\documentclass[12pt,reqno]{amsart}
\usepackage{amsmath, amssymb, amsthm}
\usepackage{hyperref}
\usepackage{enumerate}
\usepackage{mathrsfs}
\usepackage{url}
\begin{document}
\bibliographystyle{plain}
\newtheorem{thm}{Theorem}
\theoremstyle{definition}
\newtheorem{defn}[thm]{Definition}
\newtheorem{cor}[thm]{Corollary}
\newtheorem{lem}[thm]{Lemma}
\newtheorem{con}[thm]{Conjecture}
\theoremstyle{remark}
\newtheorem{rmk}[thm]{Remark}
\numberwithin{equation}{subsection}
\numberwithin{thm}{section}
%\numberwithin{defn}{section}
%\numberwithin{cor}{section}
%\numberwithin{lem}{section}
%\numberwithin{con}{section}
%\numberwithin{rmk}{section}

 \newcommand{\mo}{~\mathrm{mod}~}
 \newcommand{\mc}{\mathcal}
 \newcommand{\rar}{\rightarrow}
 \newcommand{\Rar}{\Rightarrow}
 \newcommand{\lar}{\leftarrow}
 \newcommand{\lrar}{\leftrightarrow}
 \newcommand{\Lrar}{\Leftrightarrow}
 \newcommand{\zpz}{\mathbb{Z}/p\mathbb{Z}}
 \newcommand{\mbb}{\mathbb}
 \newcommand{\B}{\mc{B}}
 \newcommand{\cc}{\mc{C}}
 \newcommand{\D}{\mc{D}}
 \newcommand{\E}{\mc{E}}
 \newcommand{\F}{\mathbb{F}}
 \newcommand{\G}{\mc{G}}
 \newcommand{\ZG}{\Z (G)}
 \newcommand{\FN}{\F_n}
 \newcommand{\I}{\mc{I}}
 \newcommand{\J}{\mc{J}}
 \newcommand{\M}{\mc{M}}
 \newcommand{\nn}{\mc{N}}
 \newcommand{\qq}{\mc{Q}}
 \newcommand{\PP}{\mc{P}}
 \newcommand{\U}{\mc{U}}
 \newcommand{\X}{\mc{X}}
 \newcommand{\Y}{\mc{Y}}
 \newcommand{\itQ}{\mc{Q}}
 \newcommand{\C}{\mathbb{C}}
 \newcommand{\R}{\mathbb{R}}
 \newcommand{\N}{\mathbb{N}}
 \newcommand{\Q}{\mathbb{Q}}
 \newcommand{\Z}{\mathbb{Z}}
 \newcommand{\A}{\mathbb{A}}
 \newcommand{\ff}{\mathfrak F}
 \newcommand{\fb}{f_{\beta}}
 \newcommand{\fg}{f_{\gamma}}
 \newcommand{\gb}{g_{\beta}}
 \newcommand{\vphi}{\varphi}
 \newcommand{\whXq}{\widehat{X}_q(0)}
 \newcommand{\Xnn}{g_{n,N}}
 \newcommand{\lf}{\left\lfloor}
 \newcommand{\rf}{\right\rfloor}
 \newcommand{\lQx}{L_Q(x)}
 \newcommand{\lQQ}{\frac{\lQx}{Q}}
 \newcommand{\rQx}{R_Q(x)}
 \newcommand{\rQQ}{\frac{\rQx}{Q}}
 \newcommand{\elQ}{\ell_Q(\alpha )}
 \newcommand{\oa}{\overline{a}}
 \newcommand{\oI}{\overline{I}}
 \newcommand{\dx}{\text{\rm d}x}
 \newcommand{\dy}{\text{\rm d}y}
 \newcommand{\cal}[1]{\mathcal{#1}}
 \newcommand{\cH}{{\cal H}}
 \newcommand{\diam}{\operatorname{diam}}

\parskip=0.5ex

\title[]{Upper bounds for double exponential sums along a subsequence}
\author{Christopher~J.~White}
\subjclass[2010]{}
\address{School of Mathematics, University of Bristol, Bristol, UK}
\email{chris.white@bristol.ac.uk}

\allowdisplaybreaks

\maketitle

\begin{abstract}
We consider a class of double exponential sums studied in a paper of Sinai and Ulcigrai. They proved a linear bound for these sums along the sequence of denominators in the continued fraction expansion of $\alpha$, provided $\alpha$ is badly-approximable. We provide a proof of a result, which includes a simple proof of their theorem, and which applies for all irrational $\alpha$.
\end{abstract}

\section{Introduction}

\subsection{Some Notation}
Let $\alpha=[a_{0};a_{1},\ldots]$ denote the continued fraction expansion of $\alpha\in\R\setminus\Z$. 
We write $||x||$ for the distance from $x$ to the nearest integer. 
The \textit{convergents} $p_{n}/q_{n}=[a_{0};a_{1},\ldots, a_{n}]$, where $(p_{n},q_{n})=1$, give good approximations to $\alpha$.
We call $\{q_{n}\}_{n\in\N}$ the \textit{sequence of denominators} of $\alpha$.
We say that an irrational number $\alpha$ is \textit{badly-approximable} if there exists $\varepsilon_{\alpha}>0$ such that for all $p,q\in\Z$, $(p,q)=1$, we have
\begin{equation}
\left|\alpha-\frac{p}{q}\right|>\frac{\varepsilon_{\alpha}}{q^{2}}.
\end{equation}
These correspond precisely with those numbers $\alpha$ for which there exists $N\in\N$ such that $a_{n}(\alpha)\leq N$ for all $n\in\N$. The set of all badly-approximable numbers is a set of Lebesgue measure zero.\\
When $\alpha$ is badly approximable, we have the helpful bound that 
\begin{equation}
||q_{n}\alpha||>\frac{\varepsilon_{\alpha}}{q_{n}}
\end{equation}
and since convergents give the best approximations for the distance to the nearest integer (see  \cite{Khin}), this means that for $m\leq q_{n+1}-1$ we have the bound
\begin{equation}
||m\alpha||>\frac{\varepsilon_{\alpha}}{q_{n}}.
\end{equation}
We write $f(n)=O(g(n))$ to mean that there exists a constant $C$ (which doesn't depend on $n$), such that $f(n)\leq C\cdot g(n)$ for all $n\in\N$.\\
Finally we define the discrepancy of a sequence.
\begin{defn}
\label{discrep}
Let $(x_{n})$ be a sequence of real numbers. For $N\in\N$ the \textit{discrepancy} of $(x_{n})$ modulo one, $D_{N}(x_{n})$, is defined as:
\begin{equation}
D_{N}(\{x_{m}\}):=\sup_{I\subseteq\mathbb{R}/\mathbb{Z}}\left|\sum_{m=1}^{N}\chi_{I}(x_{m})-N\cdot|I|\right|
\end{equation}
where $I$ denotes an interval and $\chi_{I}$ is the characteristic function of $I$.
\end{defn}

\subsection{Double Exponential Sums}
In \cite{Corinna} Sinai \& Ulcgrai studied double trigonometric sums of the form
\begin{equation}
\label{equation}
T_{M}(\alpha)=\frac{1}{M}\sum_{m=0}^{M-1}\sum_{n=0}^{M-1}{e(nm\alpha)}.
\end{equation}
We want to determine when the absolute value of this sum is bounded uniformly (i.e. by a constant which depends only on $\alpha$) over some subsequence $M\in\mathscr{A}\subseteq\N$. This will obviously depend on the Diophantine properties of $\alpha$ and the subsequence $\mathscr{A}$.\\ 
We will see that the problem of bounding this sum depends importantly on controlling sums such as
\begin{equation}
\label{simple}
\left|\frac{1}{M}\sum_{m=1}^{M}\frac{1}{\{\{m\alpha\}\}}\right|.
\end{equation}
Here 
\begin{equation}
\{\{x\}\}:=
\begin{cases}
\{x\}, & x\in [0,\frac{1}{2}],\\
\{x\}-1, & x\in (-\frac{1}{2},0),
\end{cases}
\end{equation}
is the \textit{signed fractional part} of $x\in\R$.\\
\\
In \cite{Corinna} the following is proved
\begin{thm}[Sinai, Ulcigrai]
\label{thm1}
Let $\alpha$ be badly-approximable. Consider the following double trigonometric sum:
\begin{equation}
\label{goal}
T_{M}(\alpha):=\frac{1}{M}\sum_{m=0}^{M-1}\sum_{n=0}^{M-1}e(nm\alpha).
\end{equation}
Then there exists a constant $C=C(\alpha)>0$ such that $|T_{M}|\leq C_{\alpha}$ for all $M\in\{q_{n}\}_{n\in\N}$.
\end{thm}
Our main theorem generalises this.
\begin{thm}
\label{mainthm}
Let $\alpha\in\R\setminus\Z$. Then there exists a constant $C=C(\alpha)>0$ such that 
\begin{equation}
|T_{q_{n}}|\leq C_{\alpha}\cdot\max\left\{\frac{\log (2\cdot\max_{i\leq n}\{a_{i}\})}{a_{n+1}},1\right\}.
\end{equation}
for all $n\in\N$.
\end{thm}
\begin{rmk}
By examining signs it appears that the upper bound here is close to best possible. Equation (1.13) in \cite{BHV} gives a lower bound for the largest terms in a sum that we will consider. While it's true that we use an inequality earlier in the calculation, it isn't too restrictive.
\end{rmk}

\section{Proof of main result}

\subsection{Reducing $T_{M}$}

Following the methods in \cite{Corinna} we split $T_{M}$ into two separate sums.\\
By summing the terms for $n=0,\ldots, M-1$ we can rewrite (\ref{equation}) as
\begin{equation}
T_{M}=1+\frac{1}{M}\sum_{m=1}^{M-1}{\frac{e(Mm\alpha)-1}{e(m\alpha)-1}}.
\end{equation}
Then we can write $T_{M}=1+S'_{M}-S''_{M}$ where
\begin{align}
S'_{M}:=&\frac{1}{M}\sum_{m=1}^{M-1}{\frac{e(Mm\alpha)}{e(m\alpha)-1}}\ \text{and}\label{sum1}\\
S''_{M}:=&\frac{1}{M}\sum_{m=1}^{M-1}{\frac{1}{e(m\alpha)-1}}.\label{sum2}
\end{align}
We will prove that there exist constants $C',C''\in\R$ such that 
\begin{equation}
\label{anewstep}
|S'_{q_{n}}|\leq C'\cdot\max\left\{\frac{\log (2\cdot\max_{i\leq n}\{a_{i}\})}{a_{n+1}},1\right\}
\end{equation}
and 
\begin{equation}
|S''_{q_{n}}|\leq C'',
\end{equation}
for all $n\in\N$. These constants will depend only on $\alpha$.

\subsection{The Sum $S''_{M}$ (\ref{sum2})}

Let's consider the `less intimidating' sum first.
We want to show that there exists $C''$ such that $|S''_{q_{n}}|\leq C''$ for all $n\in\N$.\\
Note that in \cite{HL}, Hardy and Littlewood prove a similar theorem.
\begin{thm}[Hardy, Littlewood]
\label{HL}
Let $\alpha$ be badly-approximable. Then there exists $C^{*}>0$ such that $|S^{''}_{M}|\leq C^{*}$ for each $M\in\N^{+}$.
\end{thm}
We proceed by calculating real and imaginary parts. 
\begin{equation}
\frac{1}{e(m \alpha)-1}=-\frac{1}{2}-\frac{i}{2} \cot(\pi m \alpha).
\end{equation}
The Taylor series expansion of $\cot x$ is
\begin{equation}
\cot x= \sum_{n=0}^{\infty}\frac{(-1)^{n}4^{n}B_{2n}}{(2n)!}x^{2n-1}=\frac{1}{x}-\frac{x}{3}+\frac{x^{3}}{45}-\ldots
\end{equation}
with radius of convergence $0<|x|<\pi$. Here $B_{n}$ is the $n$th \textit{Bernoulli number}.\\
Note that due to the symmetry of $\cot x$,
\begin{equation}\cot(\pi m\alpha)=\cot(\pi \{\{m\alpha\}\}).
\end{equation}
So we can write
\begin{equation}
\cot(\pi m\alpha)=\frac{1}{\pi\{\{m\alpha\}\}}\left(1+\sum_{n=1}^{\infty}\frac{(-1)^{n}4^{n}B_{2n}}{(2n)!}(\pi\{\{m\alpha\}\})^{2n}\right)
\end{equation}
Now the series on the right is negative and it takes values strictly between $0$ (when $\{\{m\alpha\}\} $ is close to $0$) and $-1$ (when $\{\{m\alpha\}\} $ is close to $\pm\frac{1}{2}$). 
\\
Hence in order to prove that $|S^{''}_{q_{n}}|$ is bounded by a uniform constant for all $n\in\N$, we have to prove the following
\begin{lem}
\label{new}
Let $\alpha\in\R$. Then there exists $C=C(\alpha)>0$ such that, for all $n\in \N$,
\begin{equation}
\label{stepstone}
\left|\sum_{m=1}^{q_{n}-1}\frac{1}{q_{n}\{\{m\alpha\}\}}\right|\leq C.
\end{equation}
\end{lem}
We will consider two different proofs of Lemma \ref{new}. The first is simpler, while the the latter will be applicable to $S'_{M}$ as well. The second proof is also malleable to proving Theorem \ref{HL}.

\subsection{Koksma-Hlawka Proof of Lemma \ref{new}}

Recall the Koksma-Hlawka inequality.
\begin{lem}
Let $f$ be a function with period $1$ of bounded variation. Then for every sequence $x_{m}$ and every integer $N\geq 1$, we have
\begin{equation*}
\left|\frac{1}{N}\sum_{m=1}^{N}{f(x_{m})}-\int_{0}^{1}f(x)dx\right|\leq V(f)\frac{D_{N}(x_{m})}{N},
\end{equation*}
where $V(f)$ is the total variation of the function.
\end{lem}
We wish to apply this inequality with $f(x)=\frac{1}{\{\{x\}\}}$, $x_{m}=\{m\alpha\}$ and $N=q_{n}-1$.\\
Therefore we have to restrict the domain on which we define our function, in order to ensure that it's integrable.\\ 
We're able to use the following from \cite{Rock}.
\begin{equation}
\label{useful}
\left|\alpha-\frac{p_{n-1}}{q_{n-1}}\right|>\frac{1}{2q_{n-1}q_{n}}.
\end{equation}
%\begin{rmk}
%In fact we can do better (\cite{Khin}, page $20$):
%\begin{equation}
%\left|\alpha-\frac{p_{n-1}}{q_{n-1}}\right|>\frac{1}{q_{n-1}(q_{n-1}+q_{n})}
%\end{equation}
%\end{rmk}
So for all $m\leq N=q_{n}-1$, we have 
\begin{equation}
||m\alpha||>\frac{1}{2q_{n}}.
\end{equation}
Hence we can restrict the domain of $f$ to the interval $[\frac{1}{2q_{n}},1-\frac{1}{2q_{n}}]$. Since $f$ is anti-symmetric about $1/2$, the integral above is equal to $0$.\\
The total variation, $V(f)$, of $f$ is 
\begin{equation*}
\sup_{p}\sum_{i=1}^{n_{p}}\left|\frac{1}{\{\{x_{i+1}\}\}}-\frac{1}{\{\{x_{i}\}\}}\right|,
\end{equation*}
where $\mathcal{P}$ is a partition of $[\frac{1}{2q_{n}},1-\frac{1}{2q_{n}}]$. It is maximised when we take the trivial partition (that is the two endpoints). Therefore $V(f)=4q_{n}$.\\
Finally we move on to considering the Discrepancy.\\
Lemma $5.6$ from \cite{Harman} states that
\begin{equation}
\label{Harman}
D_{N}(m\alpha)\leq 3\sum_{j=0}^{r}t_{j},
\end{equation}
where $N=\sum_{j=0}^{r}q_{r}t_{j}$. So if $N=q_{n}$ then $t_{n}=1$ and $t_{i}=0$ for all $i\neq n$. So $D_{q_{n}}(m\alpha)\leq 3$.
\\
Finally can apply all the estimates we have (with $N=q_{n}-1$ and $f$ \& $\{x_{m}\}$ as above). 
\begin{align}
\left|\sum_{m=1}^{q_{n}-1}f(x_{m})-(q_{n}-1)\int_{0}^{1}f(x)dx\right| & \leq D_{q_{n}-1}(x_{m})V(f)\\
& \leq(D_{q_{n}}(x_{m})+1)V(f)\\
& \leq 4\cdot 4q_{n}=16q_{n}
\end{align}
Hence 
\begin{equation}
\left|\frac{1}{q_{n}}\sum_{m=1}^{q_{n}-1}\frac{1}{\{\{m\alpha\}\}}\right|\leq 16.
\end{equation}
Here we used the obvious fact that $D_{M}(x_{m})\leq D_{M+1}(x_{m})+1$.

\subsection{The sum $S^{'}_{M}$ (\ref{sum1})}

We move on to considering the sum
\begin{equation}
S'_{M}:=\frac{1}{M}\sum_{m=1}^{M-1}{\frac{e(Mm\alpha)}{e(m\alpha)-1}},\label{sum3}
\end{equation}
We will write this sum as a telescoping series and take advantage of some cancellation to reduce our situation to considering the sum $S^{''}_{M}$ (\ref{sum2}).\\
\\
We write
\begin{align}
\sum_{m=1}^{M-1}\frac{e(Mm\alpha)}{e(m\alpha)-1}&=\sum_{m=1}^{M-1}(e(Mm\alpha)-e(M(m+1)\alpha))\sum_{k=1}^{m}\frac{1}{e(k\alpha)-1}\label{firstline2}\\ &+e(M^{2}\alpha)\sum_{k=1}^{M-1}\frac{1}{e(k\alpha)-1}\label{basics}.
\end{align}

We write $\alpha=\frac{p_{n}}{q_{n}}+\frac{\xi_{n}}{q_{n}q_{n+1}}$ where $\frac{1}{2}<|\xi_{n}|<1$. 
We then consider the outer part of the sum on the right hand side of (\ref{firstline2}) (for $M=q_{n}$)
\begin{align}
e(mq_{n}\alpha)-e((m+1)q_{n}\alpha) & =e(mq_{n}\psi_{n})-e((m+1)q_{n}\psi_{n})\\
& =e(mq_{n}\psi_{n})-e(mq_{n}\psi_{n})e(q_{n}\psi_{n})\\
& =(1-e(q_{n}\psi_{n}))e(mq_{n}\psi_{n}),
\end{align}
which in absolute value is less than $2\pi/q_{n+1}$.\\
Now using the triangle inequality and Lemma \ref{new} we see that (\ref{anewstep}) results from the following lemma 
\begin{lem}
\label{ost}
$\forall m\leq q_{n}-1$,
\begin{equation}
\sum_{k=1}^{m}\frac{1}{\{\{k\alpha\}\}}=O(q_{n}\max_{i\leq n}\{1,\log a_{i}\}).
\end{equation}
\end{lem}
To prove this Lemma we will need to introduce some different techniques, which will also yield a new proof of Lemma \ref{new}.

\subsection{Ostrowski Proof of Lemmas \ref{new} \& \ref{ost}}

Our alternative proof of Lemma \ref{new} will be to decompose the sum in (\ref{stepstone}) into segments where there is some obvious cancellation. 
\begin{defn}
\label{ostrowski}
Let $\alpha$ be irrational, then for every $n\in\N$ there exists a unique integer $M\geq 0$ and a unique sequence $\{c_{k+1}\}_{k=0}^{\infty}$ of integers such that $q_{M}\leq m<q_{M+1}$ and
\begin{equation}
m=\sum_{k=0}^{\infty}{c_{k+1}q_{k}},
\end{equation}
with 
\begin{align}
& 0\leq c_{1}<a_{1},\ 0\leq c_{k+1}\leq a_{k+1}\ \text{for}\ k\geq1,\\
& c_{k}=0\ \text{whenever}\ c_{k+1}=a_{k+1}\ \text{for some}\ k\geq 1,
\end{align}
and 
\begin{align}
c_{k+1}=0\ \text{for}\ k>M.
\end{align}
This is known as the \textit{Ostrowski expansion}.
\end{defn}
We will consider segments of our sum which `spread out' in the unit interval. We take our inspiration from a set of intervals discussed in \cite{Haynes}.
\begin{defn}[Special Intervals]
Define $A(m,c)$ to be the collection of non-negative integers $n$ with Ostrowski expansions of the form
\begin{equation}
n=cq_{m-1}+\sum_{k=m}^{\infty}c_{k+1}q_{k}
\end{equation}
and define a set $J(m,c)$ (which turns out be be an interval, see \cite{Haynes}) in $\R/\Z$ by
\begin{equation}
J(m,c)=\overline{\{n\alpha:n\in A(m,c)\}}
\end{equation}
For any fixed $m$, these intervals cover $\R/\Z$ and have some very nice properties such as the discrepancy of $\{n\alpha\}$ being bounded. %by a constant within them.
We will use what these intervals tell us about the distribution of $n\alpha$ on the unit interval to achieve cancellation in (\ref{simple}).
\end{defn}
Let $m\leq q_{n}-1$,
\begin{equation}
m=\sum_{i=0}^{n-1}c_{i+1}q_{i},\quad 0\leq c_{i+1}\leq a_{i+1}
\end{equation}
and 
\begin{equation}
n(i,c):=\sum_{j=0}^{i-1}c_{j+1}q_{j}+cq_{i}.
\end{equation}
We will use this decomposition to sum up to $m$.
\begin{equation}
\label{newish}
\sum_{k=1}^{m}{\frac{1}{\{\{k\alpha\}\}}}=\sum_{i=0}^{n-1}\sum_{c=0}^{c_{i+1}-1}\sum_{l=n(i,c)+1}^{n(i,c)+q_{i}}\frac{1}{\{\{l\alpha\}\}}
\end{equation}
Note that $n(i,c)+q_{i}=n(i,c+1)$ and $n(i,c_{i+1}-1)+q_{i}=n(i+1,0)$.\\
Let's consider a situation where we are studying
\begin{equation}
\label{typical}
\sum_{l=n(i,c)+1}^{n(i,c)+q_{i}}\frac{1}{\{\{l\alpha\}\}}.
\end{equation}
We wish to approximate $\alpha$ by $p_{i}/q_{i}$ and achieve (almost) complete cancellation in the main term that we get.\\
Obviously problems can occur. Specifically, if $l\cdot p_{i}\equiv 0(q_{i})$ then we don't want to divide by $0$, so we want to isolate these terms and deal with them separately. Note that since $(p_{i},q_{i})=1$, we have a complete set of residue classes modulo $q_{i}$, so in each sum, (\ref{typical}), we will have exactly one term, $(c+1)q_{i}$, where this happens .\\ 
Also, there exists $r\leq q_{i}$ such that
\begin{align}
n(i,0)+r=q_{i}\\
n(i,1)+r=2q_{i}
\end{align}
\begin{equation*}
\vdots
\end{equation*}
\begin{align}
n(i,c_{i+1}-1)+r=c_{i+1}q_{i}
\end{align}
So we can consider all these terms separately.

Finally we consider summing over a complete set of residue classes modulo $q_{i}$. We will first consider the simple case, ($1\leq k\leq q_{i-1}$), which will give us a second proof of Lemma \ref{new}.\\
We write
\begin{equation}
\alpha=\frac{p_{i}}{q_{i}}+\frac{\xi_{i}}{q_{i}q_{i+1}},
\end{equation}
where $\frac{1}{2}<|\xi_{i}|<1$. Now
\begin{align}
\label{thingy}
\sum_{k=1}^{q_{i}-1}{\frac{1}{\{\{k\alpha\}\}}}& =\sum_{k=1}^{q_{i}-1}{\frac{1}{{\{\{k\frac{p_{i}}{q_{i}}+\frac{k\xi_{i}}{q_{i}q_{i+1}}\}\}}}}.
\end{align} 
Now we use the fact that
\begin{equation}
\left\{\left\{\frac{kp_{i}}{q_{i}}+\frac{k\xi_{i}}{q_{i}q_{i+1}}\right\}\right\}=\left\{\left\{\frac{kp_{i}}{q_{i}}\right\}\right\}+\left\{\left\{\frac{k\xi_{i}}{q_{i}q_{i+1}}\right\}\right\},
\end{equation}
unless perhaps if $kp_{i}\equiv\frac{q_{i}}{2}$ modulo $q_{i}$ (when $2|q_{i}$), or if $kp_{i}\equiv\frac{q_{i}\pm 1}{2}$ (when $2|q_{i}+1$).
Now (\ref{thingy}) equals\footnote{The one or two extra term/s mentioned just above have been removed from the sum and are accounted for by the $O(1)$ term.}
\begin{align}
\sideset{}{'}\sum_{k=1}^{q_{i-1}}{\frac{1}{\{\{k\frac{p_{i}}{q_{i}}\}\}}\left(\frac{1}{1+\{\{\frac{kp_{i}}{q_{i}}\}\}^{-1}\frac{k\xi_{i}}{q_{i}q_{i+1}}}\right)}+O(1).
\end{align}
Furthermore
\begin{gather}
\left(\frac{1}{1+\{\{\frac{kp_{i}}{q_{i}}\}\}^{-1}\frac{k\xi_{i}}{q_{i}q_{i+1}}}\right)\\
=1-\left\{\left\{\frac{kp_{i}}{q_{i}}\right\}\right\}^{-1}\frac{k\xi_{i}}{q_{i}q_{i+1}}+\left\{\left\{\frac{kp_{i}}{q_{i}}\right\}\right\}^{-2}\left(\frac{k\xi_{i}}{q_{i}q_{i+1}}\right)^{2}-\ldots
\end{gather}
There exists $n_{k}$ such that $1\leq n_{k}\leq q_{i}-1$ and $n_{k}\equiv kp_{i}\mod q_{i}$. Now we define $n^{'}_{k}$ as follows
\begin{equation}
n^{'}_{k}:=
\begin{cases}
n_{k}, & n_{k}\leq \frac{q_{i}}{2}\\
n_{k}-q_{i}, & n_{k}>\frac{q_{i}}{2}.
\end{cases}
\end{equation}
Then
\begin{equation}
\left\{\left\{\frac{kp_{i}}{q_{i}}\right\}\right\}^{-1}\frac{k\xi_{i}}{q_{i}q_{i+1}}=\frac{k\xi_{i}}{n^{'}_{k}q_{i+1}}.
\end{equation}
We then know that for all $k$,
\begin{gather}
-\left\{\left\{\frac{kp_{i}}{q_{i}}\right\}\right\}^{-1}\frac{k\xi_{i}}{q_{i}q_{i+1}}+\left\{\left\{\frac{kp_{i}}{q_{i}}\right\}\right\}^{-2}\left(\frac{k\xi_{i}}{q_{i}q_{i+1}}\right)^{2}-\ldots\\
=C_{k}\frac{k\xi_{i}}{n^{'}_{k}q_{i+1}}.
\end{gather}
We need $|n^{'}_{k}|\geq 2$ in order to have a uniform bound over $k$ for the constant $C_{k}$. When this is the case 
\begin{equation}
-\frac{1}{2}<C_{k}<2,
\end{equation}
(apart from the one or two exceptions mentioned previously).
So we have to isolate another two terms.  We write $k_{1}$, $k_{-1}$ for the numbers where $k_{1}p_{i}\equiv 1\mod q_{i}$ and $k_{-1}p_{i}\equiv -1\mod q_{i}$ respectively.\\
So (\ref{thingy}) becomes
\begin{gather}
\sum_{n_{k}=2}^{q_{i}-2}{\left(\frac{1}{\{\{\frac{n_{k}}{q_{l}}\}\}}+C_{k}\left(\frac{k\xi_{i}q_{i}}{(n^{'}_{k})^{2}q_{i+1}}\right)\right)}+\frac{1}{\{\{k_{1}\alpha\}\}}+\frac{1}{\{\{k_{-1}\alpha\}\}}+O(1)\\
=\sum_{n_{k}=2}^{q_{i}-2}{C_{k}\left(\frac{k\xi_{i}q_{i}}{(n^{'}_{k})^{2}q_{i+1}}\right)}+O(q_{i}).
\end{gather}
(We used the basic approximation from Khinchin (\ref{useful}) to deal with the two extra terms.)\\
By the rearrangement inequality (see \cite{HLP}, Theorem 368) this first sum is less than 
\begin{equation}
q_{i}\left(\frac{1}{2^{2}}+\frac{1}{3^{2}}+\ldots\right),
\end{equation}
which in turn is bounded above by $q_{i}$.\\
So 
\begin{equation}
\sum_{k=1}^{q_{i}-1}{\frac{1}{\{\{k\alpha\}\}}}=O(q_{i}),
\end{equation}
as required.\\
Now we move on to a proof of Lemma \ref{ost}. We wish to prove (for all $i$) that
\begin{equation}
\sum_{c=0}^{c_{i+1}-1}\sum_{l=n(i,c)+1}^{n(i,c)+q_{i}}\frac{1}{\{\{l\alpha\}\}}=O(q_{i+1}\log c_{i+1}).
\end{equation}
Note that if we sum 
\begin{equation}
\sum_{l=n(i,c)+1}^{n(i,c)+q_{i}}\frac{1}{\{\{l\alpha\}\}}
\end{equation}
then a similar argument to the proof of Lemma \ref{new} shows that this is equal to 
\begin{equation}
O(q_{i})+\frac{1}{\{\{k_{(1,c)}\alpha\}\}}+\frac{1}{\{\{k_{(-1,c)}\alpha\}\}}+\frac{1}{\{\{(c+1)q_{i}\alpha\}\}},
\end{equation} 
where $n(i,c)+1\leq k_{(\pm 1,c)}\leq n(i,c)+q_{i}$, and $k_{(\pm 1,c)}p_{i}\equiv \pm 1\mod q_{i}$.\\
Clearly $k_{(\pm 1,c)}=k_{(\pm 1,0)}+cq_{i}$.\\
Furthermore as $m<q_{i+1}$,
\begin{equation}
k_{(\pm 1,c_{i+1}-r)}<n(i,(c_{i+1}-r))+q_{i}< q_{i+1}-(r-1)q_{i}.
\end{equation}
Now, we calculated earlier that
\begin{equation}
\frac{1}{\{\{k\alpha\}\}}=\frac{1}{\{\{k\frac{p_{i}}{q_{i}}\}\}}\left(\frac{1}{1+\{\{\frac{kp_{i}}{q_{i}}\}\}^{-1}\frac{k\xi_{i}}{q_{i}q_{i+1}}}\right)
\end{equation}
Letting $k=k_{(1,0)}$
\begin{align}
\frac{1}{\{\{k_{(1,0)}\alpha\}\}}&=q_{i}\left(\frac{1}{1+\frac{k_{(1,0)}\xi_{i}}{q_{i+1}}}\right)\\
&= \frac{q_{i}q_{i+1}}{q_{i+1}+k_{(1,0)}\xi_{i}}.
\end{align}
Hence
\begin{align}
\frac{1}{\{\{k_{(1,c)}\alpha\}\}}&=\frac{q_{i}q_{i+1}}{q_{i+1}+(k_{(1,0)}+cq_{i})\xi_{i}}
\end{align}
and also
\begin{align}
\frac{1}{\{\{k_{(-1,c)}\alpha\}\}}&=\frac{-q_{i}q_{i+1}}{q_{i+1}-(k_{(-1,0)}+cq_{i})\xi_{i}}
\end{align}
Now, without loss of generality, assume that $\xi_{i}>0$. Then 
\begin{equation}
\frac{1}{\{\{k_{(1,c)}\alpha\}\}}<q_{i}
\end{equation}
for all $c$. Hence 
\begin{align}
\sum_{c=0}^{c_{i+1}-1}\sum_{l=n(i,c)+1}^{n(i,c)+q_{i}}{\frac{1}{\{\{l\alpha\}\}}} &=O(q_{i+1})+\sum_{c=0}^{c_{i+1}-1}{\frac{1}{\{\{(c+1)q_{i}\alpha\}\}}}\\&+\sum_{c=0}^{c_{i+1}-1}{\frac{-q_{i}q_{i+1}}{q_{i+1}-(k_{(-1,0)}+cq_{i})\xi_{i}}}\\
&=O(q_{i+1})+O(q_{i+1}\log c_{i+1})\\
&+\sum_{c=0}^{c_{i+1}-1}{\frac{-q_{i}q_{i+1}}{q_{i+1}-(k_{(-1,0)}+cq_{i})\xi_{i}}}.
\end{align}
Finally 
\begin{align}
\left|\sum_{c=0}^{c_{i+1}-1}{\frac{-q_{i}q_{i+1}}{q_{i+1}-(k_{(-1,0)}+cq_{i})\xi_{i}}}\right|& \leq\frac{q_{i}}{c_{i+1}-1}+\ldots+\frac{q_{i}}{2}+q_{i}+2q_{i+1}\\
& =O(q_{i+1}\log c_{i+1}).
\end{align}
As this is true for all $i$, the condition for Lemma \ref{ost} follows.
\begin{rmk}
Equation (1.13) in \cite{BHV} tells us that the sum 
\begin{equation}
\sum_{c=0}^{c_{i+1}-1}{\frac{1}{\{\{(c+1)q_{i}\alpha\}\}}}
\end{equation}
can be no smaller than $O(q_{i+1}\log c_{i+1})$.
\end{rmk}

\begin{rmk}
In our final calculation we have ignored the cancellation between the positive and negative terms. However, when $c_{i+1}\approx a_{i+1}/2$ for example, we get very little cancellation and our main term is $O(q_{i+1}\log a_{i+1})$
\end{rmk}

\end{document}